\documentclass[12pt]{amsart}
\usepackage{geometry}
\usepackage{amsmath}
\usepackage{amsfonts}
\usepackage{amsthm}
\usepackage{amssymb}
\allowdisplaybreaks[2]

\newtheorem{theorem}{Theorem}%[section]

\newtheorem{other}{\bf Theorem}

\newenvironment{Pf}{\noindent{\emph{Proof of}}}{$\Box$}

\numberwithin{equation}{section}

\newcommand{\D}{\mathbb D}
\newcommand{\C}{\mathbb C}
\newcommand{\hol}{\mathcal Hol}
\newcommand{\ig}{\stackrel{\text{def}}{=}}

\begin{document}

\title[Bounded functions
in analytic Besov spaces]{A radial integrability result concerning
bounded functions in analytic Besov spaces with applications}

\author{Salvador Dom\'{\i}nguez}

\author{Daniel Girela}
\address{An\'{a}lisis Matem\'{a}tico, Facultad de Ciencias, Universidad de M\'{a}laga, 29071 M\'{a}laga, Spain}
\email{sdmolina16@hotmail.com} \email{girela@uma.es}
\thanks{This research is supported in part by a grant from \lq\lq El Ministerio de
Econom\'{\i}a y Competitividad\rq\rq , Spain (PGC2018-096166-B-I00)
and by grants from la Junta de Andaluc\'{\i}a (FQM-210 and
UMA18-FEDERJA-002).}

%    General info
\subjclass[2000]{30H30, 47H30; 47B38}

\keywords{Besov spaces, Radial integrability, Superposition
operators, Multipliers}

\begin{abstract} We prove that for every $p\ge 1$ there exists a
bounded function in the analytic Besov space $B^p$ whose derivative
is \lq\lq badly integrable\rq\rq\, along every radius. We apply this
result to study multipliers and weighted superposition operators
acting on the spaces $B^p$.
\end{abstract}

\maketitle

\section{Introduction}\label{intro}
Let $\D=\{z\in\C:|z|<1\}$ denote the open unit disc in the complex
plane $\mathbb C$ and let $\hol(\mathbb D)$ be the space of all
analytic functions in $\D$ endowed with the topology of uniform
convergence in compact subsets. Also, $dA$ will denote the area
measure on $\mathbb D$, normalized so that the area of $\mathbb D$
is $1$. Thus $dA(z)\,=\,\frac{1}{\pi }\,dx\,dy\,=\,\frac{1}{\pi
}\,r\,dr\,d\theta $.
\par\medskip For $0\,\le\,r\,<\,1$ and
$g$ analytic in $\mathbb D $ we set
$$M_ p(r, g)=\left (
\frac{1}{2\pi }\int_ {0}^{2\pi }\left \vert g(re\sp {i\theta })
\right \vert \sp pd\theta \right )\sp {1/p}, \quad 0<p<\infty , $$
\smallskip
$$M_ \infty (r, g)=\max_ {\vert z\vert =r}\vert g(z)\vert .$$
For $0<p\leq \infty $ the Hardy space $H\sp p$ consists of those
functions $g$, analytic in $\mathbb D $, for which
$$\left \vert \left \vert g \right \vert \right\vert _ {H\sp p}\,\ig\,\sup_ {0<r<1}
M_ p(r, g)<\infty .$$ We refer to \cite{Du:Hp} for the theory of
Hardy spaces.
\par\medskip
For $0<p<\infty $ and $\alpha >-1$ the weighted Bergman space
$A^p_\alpha $ consists of those $f\in \hol (\mathbb D)$ such that
$$\Vert f\Vert _{A^p_\alpha }\,\ig\, \left ((\alpha +1)\int_\mathbb D(1-\vert z\vert ^2)^{\alpha }\vert f
(z)\vert ^p\,dA(z)\right )^{1/p}\,<\,\infty .$$ The unweighted
Bergman space $A^p_0$ is simply denoted by $A^p$. We refer to
\cite{DS,HKZ,Zhu} for the notation and results about Bergman spaces.
\par\medskip
For $1<p<\infty $, the {\it analytic Besov space\/} $B^p $ is
defined as the set of all functions $f $ analytic in $\mathbb D $
such that $f\sp\prime \in A^p_{p-2} $. Thus a function $f\in \hol
(\mathbb D)$ belongs to $B^p$ if and only if $\rho_{_p}(f)<\infty $,
where
$$\rho_{_p}(f)\,=\,\Vert f^\prime \Vert _{A^p_{p-2}}
\,=\,\left ((p-1)\int_{\mathbb  D}(1-\vert z\vert ^2)^{p-2}\vert
f^\prime (z)\vert ^p\,dA(z)\right )^{1/p}.$$
\par
 All $B^p $ spaces ($1<p<\infty $) are
conformally invariant with respect to the semi-norm $\rho_{_p}$. An
important and well-studied case is the classical {\it Dirichlet
space\/} $B^2 $ (often denoted by $\mathcal D$) of analytic
functions whose image Riemann surface has a finite area. We mention
\cite{AFP} as a fundamental reference for Besov spaces and  Zhu's
monograph \cite{Zhu} as a very good place to find a lot of
information about them.
\par\medskip
The space $B^1 $ requires a special definition: it is the space of
all analytic functions $f $ in $\mathbb D $ for which $f^{\prime
\prime} \in A^1 $. Although the corresponding semi-norm is not
conformally invariant, the space itself is. Another possible
definition (with a conformally invariant semi-norm) is given in
\cite{AFP}, where $B^1 $ was denoted by ${\mathcal M} $. \par It is
well known that all the Besov spaces are contained in the space
$BMOA$ (even more, in $VMOA$) and, hence in the Bloch space
$\mathcal B$. We refer to \cite{G:BMOA} and \cite{ACP} for the
theory of these spaces. The inclusion $B^p\subset \mathcal B$ yields
\begin{equation}\label{incl}B^1\,\subset B^q\,\subset B^p\,\subset
\mathcal B,\quad 1\,<q\,<\,p\,<\,\infty .\end{equation}
\par\medskip
Obtaining results about the integrability along radii of distinct
classes of analytic functions in the unit disc has shown to be an
important question in complex analysis which has attracted the
attention of lots of authors over the years. One of the best known
results in this line is due to Rudin \cite{Rudin-1955} who showed
the existence of an $H^\infty $-function $f$ for which the radial
variation $V(f,\theta )=\int_0^1\vert f^\prime (re^{i\theta })\vert
dr$ is infinite for every $\theta $ except possibly for those
$\theta $ in a set of the first category and of measure zero.
Bourgain \cite{Bour} solved a question raised by Rudin by showing
that, for $f\in H^\infty $, the set of those $\theta $ for which
$V(f,\theta )$ is finite cannot be empty since it must have
Hausdorff dimension $1$.
\par\medskip In this paper we shall be concerned with radial integrability properties of $B^p$-functions. By the definition, it is clear that if $1\,<\,p\,<\,\infty \,$ and $f\in \hol (\mathbb
D)$ then
\begin{equation}\label{Bp-int-mean}f\in B^p\,\Leftrightarrow\,
\int_0^1\,(1-r)^{p-2}M_p(r, f^\prime )^p\,dr\,<\,\infty .
\end{equation} Clearly, this implies that if
$1<p<\infty $ and $f\in B^p$ then
\begin{equation}\label{Bp-rad-a.e}
\text{$\int_0^1(1-r)^{p-2}\vert f^\prime (re^{i\theta })\vert
^p\,dr\,<\,\infty $,\,\, for almost every $\theta \in [0, 2\pi
]$.}\end{equation}

\par\medskip In our first result we prove that (\ref{Bp-int-mean}) and
(\ref{Bp-rad-a.e}) are sharp in a very strong sense connecting this
with (\ref{incl}). Indeed, for $1<q<p<\infty $ we prove the
existence a of function $f\in B^p\cap H^\infty $ with
$M_p(r,f^\prime )$ \lq\lq as big as possible\rq\rq \, and having
\lq\lq bad integrability properties of order $q$ along all the
radii\rq\rq . Before stating it, let us notice that throughout the
paper we shall be using the convention that $C=C(p, \alpha ,q,\beta
, \dots )$ will denote a positive constant which depends only upon
the displayed parameters $p, \alpha , q, \beta \dots $ (which often
will be omitted) but not necessarily the same at different
occurrences. Moreover, for two real-valued functions $E_1, E_2$ we
write $E_1\lesssim E_2$, or $E_1\gtrsim E_2$, if there exists a
positive constant $C$ independent of the arguments such that
$E_1\leq C E_2$, respectively $E_1\ge C E_2$. If we have
$E_1\lesssim E_2$ and $E_1\gtrsim E_2$ simultaneously then we say
that $E_1$ and $E_2$ are equivalent and we write $E_1\asymp E_2$.

\medskip
\begin{theorem}\label{main-q>1} Suppose that $1<q<p<\infty $ and
let $\phi $ be a positive increasing function defined in $[0,1)$
satisfying
\begin{equation}\label{p-phi}\int_0^1(1-r)^{p-2}\phi
(r)^p\,dr\,<\,\infty \end{equation} and
\begin{equation}\label{q2-varphi}\int_0^1(1-r)^{q-2}\phi
(r)^q\,dr\,=\,\infty .\end{equation} Then there exists a function
$f\in B^p\cap H^\infty \setminus B^q$ with the following two
properties:
\begin{equation}\label{Mp}M_p(r,f^\prime )\,\gtrsim \phi
(r).\end{equation}
\begin{equation}\label{rad-cond-a.q}\int_0^1(1-r)^{q-2}\vert f^\prime
(re^{i\theta })\vert ^q\,dr\,=\,\infty ,\quad\text{for every $\theta
\in [0, 2\pi ]$}.\end{equation}
\end{theorem}
\par\medskip
A typical example of a function $\phi $ in the conditions of
Theorem\,\@\ref{main-q>1} is
$$\phi (r)\,=\,\frac{1}{(1-r)^\alpha },\quad 0\le r<1,$$
$\alpha $ being a real  number with $1-\frac{1}{q}\,<\,\alpha
\,<\,1-\frac{1}{p}$.
\par\medskip
The substitute of Theorem\,\@\ref{main-q>1} for $q=1$ is the
following.
\begin{theorem}\label{main-q=1} Suppose that $1<p<\infty $ and
let $\phi $ be a positive increasing function defined in $[0,1)$
satisfying (\ref{p-phi}). Then there exists a function $f\in B^p\cap
H^\infty \setminus B^1$ with the following two properties:
\begin{equation}\label{Mp}M_p(r,f^\prime )\,\gtrsim \phi
(r).\end{equation}
\begin{equation}\label{rad-cond-a.1}\int_0^1\,\vert f^{\prime \prime }(re^{i\theta })\vert \,dr\,=\,\infty ,\quad\text{for every $\theta
\in [0, 2\pi ]$}.\end{equation}
\end{theorem}
\par\medskip The proofs of Theorem\,\@\ref{main-q>1} and
Theorem\,\@\ref{main-q=1} will be presented in
Section\,\@\ref{proofs-main}. In Section\,\@\ref{appl} we shall
apply these theorems to obtain results on weighted superposition
operators and multipliers acting on the Besov spaces.

\section{Proofs of the integrability results}\label{proofs-main}
{\it Proofs of Theorem\,\@\ref{main-q>1} and
Theorem\,\@\ref{main-q=1}.} Let $\phi $ be as in
Theorem\,\@\ref{main-q>1}. Set
$$r_k\,=\,1-\frac{1}{2^k},\qquad k\,=\,1, 2, \dots .$$
Since $\phi $ is increasing it easy to see that (\ref{p-phi})
implies that
\begin{equation}\label{sum-p}\sum_{k=1}^\infty
\frac{1}{2^{k(p-1)}}\phi (r_k)^p\,<\,\infty .\end{equation} For
$k\ge 1$, define $a_k=\frac{1}{2^k}\phi (r_k)$ and set
$$f(z)\,=\,\sum_{k=1}^\infty a_kz^{2^k},\quad z\in \mathbb D.$$ Then $f$ is an analytic
function in $\mathbb D$ given by a power series with Hadamard gaps.
Using (\ref{sum-p}) we see that
$$\sum_{k=1}^\infty 2^k\vert a_k\vert ^p\,=\,\sum _{k=1}^\infty \frac{1}{2^{k(p-1)}}\phi (r_k)^p\,<\,\infty
.$$ Then it follows that $f\in B^p$ (see, e.\,\@g.,
\cite[Theorem\,\@D]{DGV1}). Let $q$ be the exponent conjugate to
$p$, that is, $q=p/(p-1)$. Using H\"{o}lder's inequality with the
exponents $p$ and $q$ and (\ref{sum-p}), we obtain
\begin{align*}\sum_{k=1}^\infty \vert a_k\vert \,=\,&
\sum_{k=1}^\infty \frac{1}{2^{k/p}2^{k/q}}\phi (r_k)\, \le \,\left
(\sum_{k=1}^\infty \frac{1}{2^{kq/p}}\right )^{1/q}\left
(\sum_{k=1}^\infty\frac{1}{2^{k(p-1)}}\phi (r_k)^p\right
)^{1/p}\,<\,\infty .
\end{align*}
Then it follows that $f\in H^\infty $ (even more, $f$ belongs to the
disc algebra).
\par Now, $zf^\prime (z)=\sum_{k=1}^\infty \phi (r_k)z^{2^k}$ ($z\in
\mathbb D)$ and then it follows that
\begin{align}\label{m2}
M_2(r,f^\prime )^2\,\ge\, \sum_{k=1}^\infty \phi
(r_k)^2r^{2^{k+1}},\quad 0<r<1.\end{align} For a given $r\in (0,1)$
take $k\in \mathbb N$ such that $r_k\le r<r_{k+1}$. Using
(\ref{m2}), the facts that the functions $M_2(\cdot ,f^\prime )$ and
$\phi $ are increasing in $(0, 1)$, and the simple estimate
$r_k^{2^k}\gtrsim 1$, we obtain
\begin{align}\label{m22}M_2(r,f^\prime )^2\,\ge \,M_2(r_k,f^\prime )^2\,\ge \,\phi
(r_{k+1})^2r_k^{2^{k+2}}\,\gtrsim \phi (r_{k+1})^2\,\ge\, \phi
(r)^2.\end{align} Now, since $f$ is given by a power series with
Hadamard gaps we have that $M_p(r,f^\prime )\asymp M_2(r,f^\prime )$
(see, e.\,\@g., Theorem\,\@8.\,\@20 in \cite[Vol.\,\@I,
p.\,\@215]{Zy}). Then (\ref{m22}) yields
\begin{align*}M_p(r,f^\prime )\,\gtrsim \phi (r).\end{align*}
Actually, we have that $M_{\lambda }(r, f^\prime)\asymp
M_2(r,f^\prime )$ for all finite $\lambda$ and, consequently, we can
assert that
\begin{align}\label{mlam}M_\lambda (r,f^\prime )\,\gtrsim \phi (r),\quad 0<\lambda <\infty .\end{align}
Using this with $\lambda =q$ and (\ref{q2-varphi}) we obtain that
$f\notin B^q$. This and \cite[Theorem\,\@D]{DGV1} imply
$\sum_{k=1}^\infty 2^k\vert a_k\vert ^q\,=\,\infty $. Then
(\ref{rad-cond-a.q}) follows using a a result of Gnuschke
\cite[Theorem\,\@1]{Gn}. This finishes the proof of
Theorem\,\@\ref{main-q>1}.
\par Since $B^1\subset B^q$, $f\notin B^1$ and, hence, $\sum_{k=1}^\infty 2^k\vert a_k\vert
\,=\,\infty $. Then the just mentioned result of Gnuschke also
yields (\ref{rad-cond-a.1}). Hence Theorem\,\@\ref{main-q=1} is also
proved.
\medskip
\section{Applications to weighted superposition operators and multipliers}\label{appl}
\par\medskip
Given an entire function $\varphi$, the superposition operator
\[S_\varphi:\hol(\D)\longrightarrow\hol(\D)\]
is defined by $S_\varphi(f)=\varphi\circ f$.
\par\medskip
The natural questions in this context are: If $X$ and $Y$ are two
normed (metric) subspaces of $\hol(\D)$, for which entire functions
$\varphi $
 does the operator $S_{\varphi}$ map $X$
into $Y$? When is $S_{\varphi }$ a bounded (or continuous) operator
from $X$ to $Y$?
\par\medskip Let us remark that we are dealing with non-linear operators.
Consequently, boundedness and continuity are not equivalent a
priori. However, Boyd and Rueda \cite{BR-14} have shown that for a
large class of Banach spaces of analytic functions $X$ and $Y$,
bounded superposition operators from $X$ to $Y$ are continuous.
\par\medskip
These questions have been studied for  different pair of spaces $(X,
Y)$. See, for example,
\cite{AMV,BV,BR-13,BR-14,BFV,CG,C,DG,GM,LZ,Ramos} and the references
therein.
\par\medskip Concerning Besov spaces a complete characterization of
the entire functions $\varphi $ which map $B^p$ into $B^q$ ($1\le
p<q<\infty $) is given in \cite{BFV}. In particular, the following
result is taken from \cite{BFV}. \begin{other}\label{BFV-BpBq}
Suppose that $1\le q<p<\infty $ and let $\varphi $ be an entire
functions. Then the superposition operator $S_\varphi $ maps $B^p$
into $B^q$ if and only if $\varphi $ is constant.
\end{other}
\par\medskip We shall use Theorem\,\@\ref{main-q>1} and
Theorem\,\@\ref{main-q=1} to obtain a new proof of
Theorem\,\@\ref{BFV-BpBq}. Actually, we shall extend this result to
the setting of weighted superposition operators.\par If $\varphi$ is
an entire function and $w\in \hol (\mathbb D)$, the {\it weighted
superposition operator}
\[S_{\varphi , w}:\hol(\D)\longrightarrow\hol(\D)\] is defined by
$$S_{\varphi , w}(f)(z)\,=\,w(z)\,\varphi \big (f(z)\big ),\quad f\in \hol
(\mathbb D),\quad z\in\mathbb D.$$ In other words, $S_{\varphi ,
w}\,=\,M_w \circ S_\varphi $, where $M_w$ is the multiplication
operator defined by $$M_w(f)(z)=w(z)f(z),\quad f\in \hol (\mathbb
D),\quad z\in\mathbb D.$$ Note that $S_\varphi =S_{\varphi , w}$
with $w(z)=1$, for all $z\in \mathbb D$. The literature on weighted
superposition operators is not as wide as the one concerning
superposition operators. We mention \cite{DG} and \cite{JWLS} as
recent papers dealing with them. We can prove the following result.

\begin{theorem}\label{weighted-super} Suppose that $1\le q<p<\infty $, $w\in \hol (\mathbb
D)$, $w\not\equiv 0$, and $\varphi $ is an entire function. Then the
weighted superposition operator $S_{w,\varphi }$ maps $B^p$ into
$B^q$ if and only if $w\in B^q$ and $\varphi $ is constant.
\end{theorem}

\par\medskip When $w\equiv 1$, Theorem\,\@\ref{weighted-super}
reduces to Theorem\,\@\ref{BFV-BpBq}.
\par\medskip
\begin{Pf}{\it Theorem\,\@\ref{weighted-super}.}
Suppose that $1\le q<p<\infty $. \par It is trivial that if $\varphi
$ is constant and $w\in B^q$ then $S_{\varphi ,w}$ maps $B^p$ into
$B^q$.
\par It is also trivial that if $\varphi $ is constant and not identically zero, and $S_{\varphi ,w}$ maps
$B^p$ into $B^q$, then $w\in B_q$.
\par\medskip
It remains to prove that if $\varphi $ is not constant, and
$S_{\varphi ,w}(B^p)\subset B^q$, then $w\equiv 0$.

Take $a\in \mathbb C$ such that $\varphi (a)\neq 0$ and let $h$ be
the constant functions defined by $h(z)=a$, for all $z\in \mathbb
D$. Since $h\in B^p$, it follows that $$S_{\varphi ,
w}(h)\,=\,\varphi (a)\cdot w\in B^q.$$ This implies that
\begin{equation}\label{wBq}w\in B^q.\end{equation}
We have to show that $w\equiv 0$. Since $B^1\subset B^s$ for all
$s>1$, it suffices to consider the case $q>1$.
\par
So, suppose that $q>1$ and  $w\not \equiv 0$. Let us use
Theorem\,\@\ref{main-q>1} to pick a function $f\in B^p\cap H^\infty
$ satisfying (\ref{rad-cond-a.q}). Since $S_{\varphi ,w}(B^p)\subset
B^q$, we deduce that $$S_{\varphi ,w}(f)\,=\,w\cdot (\varphi \circ
f)\in B^q,$$ that is,
$$\int_{\mathbb D}(1-\vert z\vert ^2)^{q-2}\left \vert
w^\prime (z)\varphi (f(z))\,+\,w(z)\varphi^\prime (f(z))f^\prime
(z)\right \vert ^q\,dA(z)\,<\,\infty .$$ Since $f\in H^\infty $, we
also have that $\varphi \circ f\in H^\infty $. This and (\ref{wBq})
imply that
$$\int_{\mathbb D}(1-\vert z\vert ^2)^{q-2}\left \vert w^\prime
(z)\varphi (f(z))\right\vert ^q\,dA(z)\,<\,\infty .$$ Then it
follows that
$$\int_{\mathbb D}(1-\vert z\vert ^2)^{q-2}\left \vert w(z)\varphi
^\prime (f(z))f^\prime (z)\right\vert ^q\,dA(z)\,<\,\infty $$ and,
hence,
\begin{equation}\label{rad-w-varprime}
\int_0^1(1-r)^{q-2}\left \vert w(re^{i\theta })\varphi ^\prime
(f(re^{i\theta }))\right \vert ^q\left \vert f^\prime (re^{i\theta
})\right \vert ^q\,dr\,<\,\infty ,\,\,\text{for almost every $\theta
\in [0, 2\pi ]$.}\end{equation} Clearly, $\varphi ^\prime \circ
f\not \equiv 0$ and  $\varphi ^\prime \circ f\in H^\infty $. These
facts and the well known inclusion $B^q\subset BMOA$  readily imply
that
\begin{equation}\label{wcotvarphiprime} w\cdot (\varphi ^\prime\circ
f)\not \equiv 0\quad \text{and \quad $w\cdot (\varphi ^\prime \circ
f)\subset H^\lambda $ for all $\lambda \in (0,\infty
)$}.\end{equation} Using Fatou's theorem and the Riesz uniqueness
theorem \cite[Chapter\,\@2]{Du:Hp}, we see that the function $w\cdot
(\varphi ^\prime \circ f)$ has a finite and non-zero radial limit at
almost every point $e^{i\theta }$ of $\partial \mathbb D$. This and
(\ref{wcotvarphiprime}) imply that
\begin{equation*}
\int_0^1(1-r)^{q-2}\left \vert f^\prime (re^{i\theta })\right \vert
^q\,dr\,<\,\infty ,\,\,\text{for almost every $\theta \in [0, 2\pi
]$.}\end{equation*} This is in contradiction with
(\ref{rad-cond-a.q}).
\end{Pf}

\par\bigskip Now we turn to apply Theorem\,\@\ref{main-q>1} to study
multipliers acting on Besov spaces. For $g\in\hol (\D )$, the
multiplication operator $M_g$ is defined by \[ M_g(f)(z)\ig
g(z)f(z),\quad f\in \hol (\D),\,\, z\in \D.
\]
If $X$ and $Y$ are two spaces of analytic function in $\mathbb D $
(which will always be assumed to be Banach or $F$-spaces
continuously embedded in $\hol (\mathbb D )$) and $g\in\hol (\mathbb
D )$ then $g$ is said to be a {\bf multiplier} from $X$ to $Y$ if
$M_g(X)\subset Y$. The space of all multipliers from $X$ to $Y$ will
be denoted by $M(X,Y)$ and $M(X)$ will stand for $M(X,X)$. Using the
closed graph theorem we see that if $g\in M(X, Y)$ then $M_g$ is a
bounded operator from $X$ into $Y$.
\par\medskip It is well known
that if $X$ is nontrivial then $M(X)\subset H^\infty $ (see,
e.\,\@g., \cite[Lemma\,\@1.\,\@1]{ADMV} or
\cite[Lemma\,\@1.\,\@10]{Vi}). Clearly, this implies the following:
\begin{equation}\label{M(X,Y)Hinfty} \text{If $Y$ is nontrivial and $Y\subset
X$ then $M(X, Y)\subset H^\infty $.}\end{equation}
\par\medskip The spaces of multipliers $M(B^p, B^q)$ have been studied in a
good number of papers. The following result is part of
\cite[Theorem\,\@2]{GaGiPe-TAMS-2011}.
\begin{other}\label{MBpBq} If $1\,\le\, q\,<\,p\,<\,\infty $\, then $M(B^p, B^q)=\{
0\}$.
\end{other}
\par The proof of this result in \cite{GaGiPe-TAMS-2011} uses, among
other facts, a decomposition theorem for Besov spaces and
Khinchine's inequality. We will show next that
Theorem\,\@\ref{MBpBq} can be obtained as a consequence of
Theorem\,\@\ref{main-q>1}.
\par\medskip
\begin{Pf}{\it Theorem\,\@\ref{MBpBq}.}
Since $B^1\subset B^s$ for all $s>1$, it suffices to prove the
result in the case $1<q<p<\infty $. So, assume this and that
$M_g(B^p)\subset B^q$.
\par Suppose that $g\not\equiv 0$.
\par Since the constant function $1$ belong to $B^p$, we see that $g\in B^q$. Also, (\ref{M(X,Y)Hinfty}) and the inclusion $B^q\subset B^q$
imply that $g\in H^\infty $. Thus we have
\begin{equation}\label{gBqHinfty}g\in B^q\cap H^\infty
.\end{equation}
 Use
Theorem\,\@\ref{main-q>1} to pick a function $f\in B^p\cap H^\infty
$ satisfying (\ref{rad-cond-a.q}). Since $M_g(B^p)\subset B^q$, we
have that $M_g(f)\,=\,g\cdot f\,\in B^q$, that is
\begin{equation}\label{fgBq}\int_{\mathbb D}(1-\vert
z\vert^2)^{q-2}\vert g^\prime (z)f(z)\,+\,g(z)f^\prime (z)\vert
^q\,dA(z)\,<\,\infty .\end{equation} Since $g\in B^q$ and $f\in
H^\infty $ we see that
$$\int_{\mathbb D}(1-\vert
z\vert^2)^{q-2}\vert g^\prime (z)f(z)\vert^q\,dA(z)\,<\,\infty .$$
This and (\ref{fgBq}) imply that
\begin{equation*}\int_{\mathbb D}(1-\vert
z\vert^2)^{q-2}\vert g(z)f^\prime (z)\vert ^q\,dA(z)\,<\,\infty
.\end{equation*} and, hence,
\begin{equation}\label{rad-q-Mg}\int_0^1(1-r)^{q-2}\vert
g(re^{i\theta })f^\prime (re^{i\theta })\vert ^q\,d\theta
\,<\,\infty ,\quad\text{for almost every $\theta \in [0, 2\pi
]$.}\end{equation} Since $g\in H^\infty $ and $g\not\equiv 0$, $g$
has a finite non-zero radial limit at almost every point $e^{i\theta
}$ of $\mathbb D$. This and (\ref{rad-q-Mg}) imply that
\begin{equation*}\int_0^1(1-r)^{q-2}\vert
f^\prime (re^{i\theta })\vert ^q\,d\theta \,<\,\infty
,\quad\text{for almost every $\theta \in [0, 2\pi
]$.}\end{equation*} This is in contradiction with
(\ref{rad-cond-a.q}).
\end{Pf}
\par\medskip {\bf Acknowledgements.} We wish to thank the referee for
reading carefully the paper and making a number of nice suggestions
to improve it.
\medskip
\bibliographystyle{amsplain}

\end{document}